\documentclass[ reprint,
superscriptaddress,
amsmath,amssymb,
aps,
nofootinbib]{revtex4-2}

\usepackage{graphicx}
\usepackage{dcolumn}
\usepackage{bm}

\usepackage{xcolor}
\usepackage{mathrsfs}
\usepackage{mathtools}

\begin{document}
\preprint{APS/123-QED}

\title{Manifold Solutions to Navier--Stokes Equations}

\author{David V. Svintradze}

\email[Correspondence email address: ]{dsvintradze@newvision.ge}

\affiliation{School of Medicine, New Vision University, Bokhua 11, 0159 Tbilisi, Georgia}

\affiliation{Niels Bohr Institute, University of Copenhagen, Blegdamsvej 17, 2100 Copenhagen, Denmark}



\date{\today}

\begin{abstract}

We have developed dynamic manifold solutions for the Navier-Stokes equations using an extension of differential geometry called the calculus for moving surfaces. Specifically, we have shown that the geometric solutions to the Navier-Stokes equations can take the form of fluctuating spheres, constant mean curvature surfaces, generic wave equations for compressible systems, and arbitrarily curved shapes for incompressible systems in various scenarios. These solutions apply to predominantly incompressible and compressible systems for the equations in any dimension, while the remaining cases are yet to be solved. We have demonstrated that for incompressible Navier-Stokes equations, geometric solutions are always bound by the curvature tensor of the closed smooth manifold for every smooth velocity field. As a result, solutions always converge for systems with constant volumes.
\end{abstract}

\maketitle

\textit{Introduction--} The Navier-Stokes (NS) equations are well-known and widely studied in fluid dynamics. However, solving these equations for generic cases is a severe problem. To address this problem, we use an extension of differential geometry to moving manifolds known as the calculus for moving surfaces (CMS). This approach allows us to solve the equations for extreme limiting cases. We focus on a moving fluid volume enclosed by a moving surface, which acts as an interface separating the fluid from the environment. By projecting the motion of the fluid onto the surface, the system gains a moving shape. The volume and the surface are treated as a single system, and finding a solution to the NS equation determines the shapes in motion or vice versa. 

More specifically, let us assume a moving surface surrounds the fluid object, and the whole system is subject to inertial motion. Then, the system's action is: 
\begin{align}
\mathcal{L}&=\mathcal{L}_S+\mathcal{L}_\Omega \nonumber \\
\mathcal{L}_S&=\int_S\frac{\rho_S V_S^2}{2}dS-\int_\Omega P_S d\Omega \label{surface action} \\
\mathcal{L}_\Omega&=\int_\Omega\frac{\rho_\Omega V_\Omega^2}{2}d\Omega-\int_\Omega E_\Omega d\Omega \label{inertia action}
\end{align}
where $\rho_S$ represents the surface mass density, $\rho_\Omega$ denotes the volumetric mass density, $V_S$ stands for the surface velocities, $V_\Omega$ indicates the system's inertial velocity in a specific inertial reference frame, $E_\Omega$ is the energy defining inertial motion, and $P_S$ refers to the volumetric energy (surface pressure) inducing surface motion. Here, $\Omega$ refers to the spatial system, while $S$ represents the surface boundary throughout this paper. Following the minimum action principle, the variation of the system action must vanish, therefore:
\begin{equation}
\delta\mathcal{L}_S=-\delta\mathcal{L}_\Omega \label{action}
\end{equation}
The variation of the $\mathcal{L}_\Omega$ inertial term leads to the NS equations, while the variation of the $\mathcal{L}_S$ surface action gives uncoupled moving manifold equations. Therefore, the solution to the $\delta\mathcal{L}_\Omega=0$ NS must be the $\delta\mathcal{L}_S=0$ moving manifold's equations solutions and vice versa. In this paper, we first clarify that variation of the $\delta\mathcal{L}_\Omega$ inertial term indeed leads to the NS equations. Then we, as a reminder, re-clarify that variation of $\delta\mathcal{L}_S$ dynamic manifolds leads to our surface dynamics equations \cite{Svintradze2017, Svintradze2018}. Next, we demonstrate the specific solutions to the equations of the moving manifold we derived earlier \cite{Svintradze2019, Svintradze2020, Svintradze2023} that also must be the NS solutions. Namely, we show that approximate solutions to the equations are fluctuating spheres along with generic wave equations when the surface normal velocity $C$ largely dominates over surface tangent velocities $V_i$. The exact solutions become practically any shape if the system is incompressible or tangent velocities largely dominate on normal ones. This implies that the solutions to incompressible Navier-Stokes equations are always bounded above by the curvature tensor of a moving, smooth, closed manifold. 

Considering that the reader might need to become more familiar with mathematical formalism, we re-introduce the basics of the extension of differential geometry to moving manifolds, known as CMS, in the appendix. We did not directly contribute to the development of mathematical formalism. Our role is in application to physical problems \cite{Svintradze2017, Svintradze2018, Svintradze2019, Svintradze2020, Svintradze2023}, see also conference abstracts that preceded publication of papers \cite{Svintradze2024a, *Svintradze2023a, *Svintradze2022a, *Svintradze2021a, *Svintradze2020a, *Svintradze2017a, *Svintradze2016a, *Svintradze2015a, *Svintradze2014a, *Svintradze2013a, *Svintradze2011a, *Svintradze2009a}. Next, we proceed with providing narratives of taking variations of (\ref{surface action}, \ref{inertia action}), even though for (\ref{inertia action}) is largely known, and for (\ref{surface action}) is published \cite{Svintradze2017, Svintradze2018}. Then, we provide specific solutions to the $\delta\mathcal{L}_S=0$ surface dynamics equations and, therefore, according to (\ref{action}), show some specific solutions to the NS problem. The solutions presented here partially and approximately solve NS for dominantly compressible systems while exactly solving the NS problem for incompressible systems.

The paper is written using tensorial notations and follows the Einstein summation convention. The appearance of the same upper and lower indices indicates summation. Greek indices generally indicate ambient space tensorial components, while Latin indices indicate surface tensors. Capital letters in lower indices stand for designations and are unrelated to tensorial nomenclature.

\textit{Calculus for Moving Surfaces (CMS) --} To briefly define the extension of differential geometry, called the calculus for moving surfaces, we note that the extension was initiated by Hadamard nearly a century ago and developed by a generation of mathematicians \cite{Hadamard, *Thomas, *MGrinfeld, Grinfeld2013}. Following the CMS developments, P. Grinfeld proposed equations of motions for massive thin fluid films \cite{Grinfeld2013, Grinfeld2010, *PhysRevLett.105.137802, *Grinfeld2012} (see Appendix A, B).  We have generalized P. Grinfeld equations for any manifolds by deriving surface dynamics equations (referred to as moving manifolds equations as well) \cite{Svintradze2017, Svintradze2018}, solved shape-dynamics problems for water drops \cite{Svintradze2019}, and extended Young-Laplace, Kelvin, and Gibbs-Thomson equations to arbitrarily curved dynamic surfaces \cite{Svintradze2020, Svintradze2023}, see also meeting abstracts that preceded the papers \cite{Svintradze2023a, *Svintradze2022a, *Svintradze2021a, *Svintradze2020a, *Svintradze2017a, *Svintradze2016a, *Svintradze2015a, *Svintradze2014a, *Svintradze2013a, *Svintradze2011a, *Svintradze2009a}. Here, we claim that those solutions are also some particular and generic solutions to the NS problem.

\textit{Shape Dynamics Equations--} One can derive surface dynamics equations based on the definitions of velocity field and spatiotemporal covariant derivatives, along with integrations theorems given in CMS textbook \cite{Grinfeld2013}. Derivation requires all basic knowledge of differential geometry and its CMS extension. To avoid repetition of algebraic manipulations but maintain the consistency of the paper, we are giving only short narratives here. For detailed derivation, see \cite{Svintradze2017, Svintradze2018} and for extensions, see \cite{Svintradze2019, Svintradze2020, Svintradze2023}. 

Let us start with the most generic form of the (\ref{surface action}) action
\begin{equation}
\mathcal{L}_S=\int_S\frac{\rho V_S^2}{2}dS-\int_S\Lambda_S dS-\int_\Omega P_S d\Omega \label{boundary action}
\end{equation}    
where $\Lambda_S$ is the surface energy density and is explicit here $\Lambda_S\notin P_S$, while is not explicit in (\ref{surface action}) $\Lambda_S\in P_S$.

Next,  we take variations of the kinetic and energy terms from surface action (\ref{boundary action}). 
\begin{equation}
\delta\mathcal{L}_S=\delta\int_S\frac{\rho_S V_S^2}{2}dS-\delta(\int_\mathbb{S}\Lambda_S dS+\int_\Omega P_S d\Omega)
\end{equation}
Detailed derivations are given in  \cite{Svintradze2017, Svintradze2018} and in its current form published in \cite{Svintradze2020, Svintradze2023}. According to integration theorems (\ref{surface integral}, \ref{space integral}) and calculations of the covariant time derivative (\ref{covariant t}) \cite{Svintradze2017, Svintradze2018, Svintradze2020, Svintradze2023}, we have  
\begin{align}
&\delta\int_S\frac{\rho_S V_S^2}{2}dS=\int_S \rho_S C(\dot\nabla C+2V^i\nabla_i C+V^a V^b B_{ab}) \nonumber\\ 
&+\rho_S V_a(\dot\nabla V^a+V^k\nabla_k V^a-C\nabla^a C-CV^k B_k^a)dS \label{total variation} \\
&\delta\int_S\Lambda_S dS=\int_S(\partial_t\Lambda_S-V^k\nabla_k\Lambda_S-\Lambda_S CB_i^i)dS \label{tension variation} \\
&\delta\int_\Omega P_S d\Omega=\int_\Omega \partial_tP_Sd\Omega+\int_SP_SCdS \label{pressure variation}
\end{align}
Here, $\rho_S$ is the surface mass density, $C$ is the surface normal velocity, $V^i$ is the surface tangent velocities, $V^\alpha$ is the ambient $\alpha=1,2,3...$ component of the $V_S$ surface velocity,  $\nabla_i,\dot\nabla$ stands for tangent space and time covariant derivatives, and small Latin indexes $i=1,2...$ refer to the surface-associated tensors. Note that the first integral in (\ref{total variation}) represents the normal component of the variation, while the second integral represents the tangent component.
$\partial_t\Lambda_S$ is mainly governed by normal deformations while $V^k\nabla_k\Lambda_S$ is tangential one so that $N^\alpha\partial_t\Lambda_S=\partial_tF^\alpha$, where $F_S^\alpha$ is some force acting on the surface in $\alpha=1,2,3$ spacial direction. A variation of the (\ref{pressure variation}) is normal because the surface pressure acts in the normal direction. 

Because (\ref{tension variation}) has mixed normal and tangent components and (\ref{pressure variation}) has only normal variation, by applying the Gauss theorem to the normal components and taking into account the conservation of surface mass the surface dynamics equations $\delta\mathcal{L}_S=0$ follow:
\begin{align}
&\dot\nabla\rho_S+\nabla_i(\rho_S V^i)=\rho_S CB_i^i \nonumber \\
&\partial_\alpha[V^\alpha(\rho_S(\dot\nabla C+2V^i\nabla_i C+V^a V^b B_{ab}) -P_S+\Lambda_S B_i^i)] \nonumber \\
&=\partial_t P_S+\partial_\alpha\partial_tF_S^\alpha \label{MME} \\
&\rho_S(\dot\nabla V^i+V^k\nabla_k V^i-C\nabla^i C-CV^k B_k^i)=-\nabla^i\Lambda_S \nonumber
\end{align}
The first equation expresses the surface mass conservation, the second depicts manifold evolution in the normal direction, and the last illustrates manifold dynamics in the tangent space.  Since $i$ is the manifold's dimension, the number of the last equation also indicates the dimensions of the tangent space. In this framework, the manifold is one less-dimensional surface than the space, and the number of dimensions is not limited, even though here we discuss the simplest case, a two-dimensional surface embedded in three-dimensional space. The analyses and results can be readily extended to any dimensional manifold. Details of derivations are given in Appendix C. 

\textit{Navier-Stokes (NS) Equations--} In the previous section, we presented the derivation narratives of the surface dynamics equations. This section derives the NS equations from $\delta\mathcal{L}_\Omega=0$, even though the derivation of the NS equations has been well established and does not require any special attention \cite{landau2013fluid}. 
Considering equation (\ref{inertia action}) and the continuity equation for mass density, we have
\begin{align}
\partial_t\rho_\Omega&=-\partial_\alpha(\rho_\Omega V_\Omega^\alpha) \nonumber \\
\delta\int_\Omega \frac{\rho_\Omega  V^2_\Omega }{2}d\Omega &=\int_\Omega\rho_\Omega \bm{V_\Omega}( V_\Omega^\alpha\partial_\alpha \bm{V_\Omega}+\partial_t \bm{V_\Omega})d\Omega \label{kinetic term}
\end{align}
See derivation in Appendix D. This is exactly the kinetic term in the NS equation. The surface ambient velocities for inertial motion $V_\Omega^\alpha$ is the $\alpha$ component of $\bm{V_\Omega}$ in flat ambient space, and the surface normal velocity for inertial motion is $\bm{V_N}=\bm{V_\Omega}\cdot\bm{N}$. Note here that there is no immediate connection between $\bm{V_\Omega}$ and $\bm{V_S}$ because in the laboratory reference frame $\bm{V_\Omega}=0$, while if the surface is in motion, then there is no reference frame in which $\bm{V_S}$ is zero.

Next, we take a variation of the second part of the (\ref{inertia action}) integral. Here, we follow Landau\&Lifshitz designation because this part does not need any special physical explanation and can be readily obtained by modeling inertial force \cite{landau2013fluid}. Let's assume there is no energy in/out flux of the system. That is, the system is "passive" fluid
\begin{equation}
\int_SE_\Omega V_NdS=\int_SE_\Omega \bm{V_\Omega}\cdot\bm{N}dS=\int_\Omega\nabla_\alpha (EV_\Omega^\alpha) d\Omega=0 \label{energy flux}
\end{equation}
where $\bm{J}=E\bm{V_\Omega}=0$ is the energy flux. We applied the Gauss theorem to convert the surface integral to the space integral. Taking the (\ref{energy flux}) and the integration theorem (\ref{space integral}) into account, we have
\begin{equation}
\delta_t\int_\Omega E_\Omega d\Omega=\int_\Omega \partial_tE_\Omega d\Omega=\int_\Omega V_\Omega^\alpha\partial_\alpha E_\Omega d\Omega \label{force}
\end{equation}
Because the negative gradient of the energy is the force, we can model it by the momentum flux density tensor $\partial_\alpha E_\Omega=-\partial_\beta M_\alpha^\beta$ where
\begin{equation}
M_\alpha^\beta=p\delta_\alpha^\beta+\rho_\Omega V_\alpha V^\beta-\sigma_\alpha^{\beta\prime} \label{momentum flux}
\end{equation}
$\sigma_{\alpha\beta}=-p\delta_{\alpha\beta}+\sigma_\alpha^{\beta\prime}$ is called stress tensor and $\sigma_\alpha^{\beta\prime}$ is the viscous stress tensor, $\delta_\alpha^\beta$ is the Kronecker delta. The equation (\ref{momentum flux}) demonstrates that the potential field's pressure difference directly affects the momentum flux density. This pressure difference triggers inertial motion, which relies on the motion's speed in the $\alpha$ and $\beta$ directions and the internal frictional force. By considering angular velocities and appropriately modeling the viscous stress tensor \cite{landau2013fluid}, one can derive (the derivation is given in the appendix) the following:
\begin{align}
V^\alpha \partial_\alpha E_\Omega&=\bm{V_\Omega}\bm{\nabla}E_\Omega \label{energy variation} \\
&=\bm{V_\Omega}\cdot(-\bm{\nabla}p+\mu\Delta \bm{V_\Omega}+(\xi+\frac{1}{3}\mu)\bm{\nabla}\partial_\alpha V_\Omega^\alpha)  \nonumber
\end{align}
(\ref{energy variation}) in combination with (\ref{kinetic term}) gives a variation of the (\ref{inertia action}) action, and forcing it to vanish trivially gives the NS equation with consideration that integral equations stand for all inertial velocity $\bm{V_\Omega}$
\begin{align}
\rho_\Omega &(\partial_t \bm{V_\Omega}+V_\Omega^\alpha\partial_\alpha \bm{V_\Omega})=\nonumber \\
&-\bm{\nabla}p+\mu\Delta \bm{V_\Omega}+(\xi+\frac{1}{3}\mu)\bm{\nabla}\partial_\alpha V_\Omega^\alpha \label{NS}
\end{align}

\textit{Solutions--} Based on the identity established above (\ref{action}), the solution to the NS equation is also the solution to the surface dynamics equations and vice versa. This means we can provide simple solutions to the surface dynamics equation for two extreme cases, which will also be NS solutions. In the first case, when the normal velocity $C$ is much larger than the tangent velocity $V^i$, we refer to it as a predominantly compressible system. In the second case, when the tangent velocity $V^i$ is much larger than the normal velocity $C$, it is a predominantly incompressible system. 

\textit{Constant Field Approximations and Dominantly Compressible Systems $C>>V^i$--} Here, we show the solution to (\ref{MME}) therefore to NS in the special case when $P_S, \Lambda_S, \bm{F}_S=const$. Then, equations (\ref{MME}) are simplified to the
\begin{align}
&\dot\nabla\rho_S+\nabla_i(\rho_S V^i)=\rho_S CB_i^i \nonumber \\
&\rho_S(\dot\nabla C+2V^i\nabla_i C+V^a V^b B_{ab})=-\Lambda_S B_i^i\label{GE} \\
&\rho_S(\dot\nabla V^i+V^k\nabla_k V^i-C\nabla^i C-CV^k B_k^i)=0 \nonumber
\end{align}
(\ref{GE}) is generally known as Grinfeld equations for thin fluid films \cite{Grinfeld2013}. Now let's assume that Grinfeld equations are approximate solutions of surface dynamics equations (\ref{MME}) for quasi-constant $P_S, \Lambda_S, \bm{F}_S\neq const$ fields. Then for homogenous surfaces $\nabla^i\Lambda_S=0$ and dominantly compressible cases $C>>V^i$, when tangential velocities are infinitesimally small $V^i\approx 0$,  the dynamic manifolds equations (\ref{MME}) simplify as 
\begin{align}
\partial_t\rho_S&=\rho_SCB_i^i\nonumber \\
-\partial_\alpha(P_SV^\alpha)&=\partial_t P_S+\partial_\alpha\partial_tF_S^\alpha \label{solved} \\
\rho_S\partial_tC&=-\Lambda_SB_i^i\nonumber
\end{align}   
Note here that for homogeneous surfaces ${\partial_\alpha P_S, \partial_\alpha F_S^\alpha=0}$ and the immediate solution to the second equation, therefore to the system (\ref{solved}) is
\begin{align}
\partial_\alpha V^\alpha&=0 \label{sphere} \\
\partial^2_tC&=-\gamma (\Delta C-2KC)\label{wave}
\end{align}
(\ref{wave}) provides the wave equation for near-planar surfaces where Gaussian curvature $K\approx 0$ is infinitesimally small. It is noteworthy that, in this particular case, even though the system is predominantly compressible, it has an approximate solution (\ref{sphere}) characteristic of dominantly incompressible systems. The first equation gives the spherical solution of
\begin{equation} 
\bm{V}=k\bm{R}/R^3 \label{spherical solution}
\end{equation} 
$k$ is some coupling constant, and $\bm{R}$ is the position vector  (see derivations in Appendix E). Taking (\ref{spherical solution}) into account, one can extract an analytical solution to (\ref{sphere}) having the form of
\begin{align}
\bm{R}=A(\bm{R}_0+\omega_\alpha R_\alpha \bm{S}^\alpha t)+BR_{0\alpha}e^{\omega_\alpha \bm{S}^\alpha t}+\psi_\alpha\bm{S}^\alpha \label{SF}
\end{align}
where $R_{0\alpha}$ is arbitrary coordinates of the initial position vector $\bm{R}_0$, $\omega_\alpha$ is the coupling constant indicating the frequency of oscillation in $\alpha$ directions, and $\bm{S}^\alpha$ is the vector component of $\bm{S}$ unit vector, $A,B$ are some constants that can be defined by initial conditions, and ${{\boldsymbol{\psi }}}_{\xi }(\xi ,t)$ is the wave function defined in the Appendix E \cite{Svintradze2019}. All derivations are in the appendix. 


\textit{Dominantly Incompressible Systems, Tangent Velocity Domination $C<<V^i$ --} Let's assume another extreme scenario where tangent velocities $V^i$ largely dominate over the normal velocity $C\approx 0$. In these predominantly incompressible systems, the surface dynamics equations (\ref{MME}) significantly simplify
\begin{align}
&\dot\nabla\rho_S+\nabla_i(\rho_S V^i)=0 \nonumber \\
&\partial_\alpha[V^\alpha(\rho_SV^a V^b B_{ab}-P_S+\Lambda_S B_i^i)]=\partial_t P_S+\partial_\alpha\partial_tF_S^\alpha \nonumber \\
&\rho_S(\dot\nabla V^i+V^k\nabla_k V^i)=-\nabla^i\Lambda_S \label{MMEI}
\end{align}
The first and the last equations restrict tangent velocities (\ref{MMEI}), while the second equation lands the analytic solution of
\begin{equation}
B_{ab}=(2P_S+\partial_\alpha F_S^\alpha-\Lambda_SB_i^i)(\rho_SV_aV_b)^{-1}\label{isolution}
\end{equation} 
We refer to the solution as the generalized Young-Laplace law for nearly incompressible cases \cite{Svintradze2023}. If one assumes an infinitesimally thin surface with infinitesimally small non-zero surface mass density, then the first and last equations from (\ref{MMEI}) become identities for homogeneous surfaces $-\nabla^i\Lambda_S=0$, therefore relaxing constraints to tangent velocities and any moving surface with (\ref{isolution}) curvature tensor for $\forall V^i$ becomes the estimated solution. 

The solution (\ref{isolution}) extends to the generic Kelvin equation when the system shape is defined by vapor $p_v$ and saturation $p_s$ pressures, then
\begin{equation}
B_{ab}=(\frac{k_BT}{v_m} \ln \frac{p_v}{p_s} +\partial_\alpha F_S^\alpha-\Lambda_SB_i^i)(\rho_S V_aV_b)^{-1} \label{Kelvin}
\end{equation}
$p_v/p_s$ is relative humidity, and $v_m$ is the molar volume of the surface (see derivation in \cite{Svintradze2020}). Note that homogenous and static shapes with spherical geometry give the Kelvin equation \cite{Svintradze2020, Svintradze2023}.

Analytically, solutions to the Navier-Stokes equations can be derived if the temperature distribution is provided instead of surface pressure distributions. In this case, the solution takes the form of generic Gibbs-Thomson equations \cite{Svintradze2023}.
\begin{equation}
B_{ab}=(\gamma H+v_m(\partial_\alpha F_S^\alpha-\Lambda_S B_i^i))(\mu V_aV_b)^{-1} \label{Gibbs-Thomson}
\end{equation}
where $\gamma=1-T/T_0$, $T$ is the temperature distribution on the surface, $T_0$ is the bulk temperature, $\mu=\rho_Sv_m$ is the surface molar mass ($v_m$ is the molar volume of the surface) and $H$ is fusion enthalpy. If $\gamma, H$ are given, one can find the surface curvature tensor and identify the moving manifold solution to NS. 

(\ref{isolution}, \ref{Kelvin}, \ref{Gibbs-Thomson}) solutions indicate that NS equations for incompressible fluids are always bound from above by the curvature tensor of the smooth closed manifolds. Therefore, the NS solution converges for $\forall \bm{V_\Omega}$ inertial velocities.

\textit{Conclusion--} We have solved the Navier-Stokes equations for two extreme cases. Namely, we have demonstrated that for predominantly compressible systems, the estimated solutions are generic wave equations and fluctuating spheres. On the other hand, the formation of any moving shapes with any surface patterns is possible for predominantly incompressible fluids. We have presented solutions for the cases where $C >> V^i$ (dominantly compressible) and $C << V^i$ (dominantly incompressible). Despite imposing restrictions on the tangent velocities (\ref{MMEI}), the equations (\ref{isolution}, \ref{Kelvin}, \ref{Gibbs-Thomson}) clarify that closed and smooth manifolds always bound the solutions of the NS equations. Therefore, we refer to them as manifold solutions. However, the case for all $\forall C, V^i$ remains unsolved mainly, even though we have fully solved the NS problem for incompressible fluids and produced partial, approximated solutions for compressible ones.

\begin{acknowledgements}
This work is supported by Shota Rustaveli Science Foundation of Georgia by grant no. STEM-22-365 and Euro Commission's ERC Erasmus+ collaborative linkage grant between the University of Copenhagen and New Vision University. We also benefit from New Vision University’s internal funding. We want to acknowledge the warm hospitality of the Biocomplexity Department of Niels Bohr Institute, University of Copenhagen, where the paper was initiated and completed. 
\end{acknowledgements}   

\appendix

\section{Gaussian Differential Geometry}
Here, we revive the concept of extrinsic differential geometry, called Gaussian geometry. Gauss's approach, in contrast to Riemannian geometry, is fundamentally extrinsic. 
Embedding in ambient space is necessary to understand better moving manifolds. Therefore, instead of describing Riemannian geometry, we start from the Gaussian one, considering that the formalism is sufficient for higher than two-dimensional surfaces. In mathematics, higher dimensional surfaces are called hyper-surfaces. To minimize terminology, we refer to any dimensional manifold as a surface and explicitly indicate whether the referred surface is two-dimensional. 

Let $S^i\in\mathbb{S}^n, i=0,1,2,...,n,X^\alpha\in\mathbb{R}^{n+1}, \alpha=0,1,2...,n+1$ be the generic coordinates of the point on the hyper-surface $\mathbb{S}$ (referred as surface) embedded on one higher dimensional flat space $\mathbb{R}$ so that $\bm{R}$ is a position vector $\bm{R}\in\mathbb{R}^{(n+1)}, \bm{R}=\bm{R}(X^\alpha, t)=\bm{R}(S^i, t)$. Greek indexes stand for the space, and Latin for the surface. Bold letters indicate vectors. Partial derivatives of the position vector define covariant base vectors of ambient space and embedded surface, while the scalar product of base vectors are metric tensors
\begin{equation}
\begin{matrix*}[l] 
\bm{X}_\alpha=\partial_\alpha\bm{R}, & \bm{S}_i=\partial_i\bm{R} \\
X_{\alpha\beta}=\bm{X}_\alpha\cdot \bm{X}_\beta=\delta_{\alpha\beta}, & S_{ij}=\bm{S}_i\cdot \bm{S}_j
\end{matrix*} \label{definition}
\end{equation}
where $\partial_\alpha=\partial/\partial X^\alpha, \partial_i=\partial/\partial S^i$. Note that since ambient space is flat $X_{\alpha\beta}=\delta_{\alpha\beta}$. Inverse to surface $S_{ij}$ tensors are contravariant metric tensors $S_{ij}S^{jk}=\delta_i^k$. The surface normal $\bm{N}$ is defined as the normal vector to surface bases so that $\bm{N}\cdot \bm{S}_i=0, \bm{N}^2=1$ and to ensure convex geometry being defined as positive curvature, the direction of the normal is set toward the surface. 

The definition of base vectors allows the introduction of Christoffel symbols for embedded surfaces. For ambient flat space, the symbols $\Gamma_{\alpha\beta}^\gamma=\bm{X}^\gamma\cdot\partial_\alpha\bm{X}_\beta=0, \partial_\alpha=\nabla_\alpha$ are vanishing, while for embedded surface $\Gamma_{ij}^k=\bm{S}^k\cdot\partial_i\bm{S}_j\neq 0$ and $\partial_i\neq\nabla_i$. 

While all symbols are unchanged in differential geometry for moving manifolds, they are not the only ones and are enhanced by extrinsic parameters. Christoffel symbols allow the definition of invariant to reference frame covariant derivatives
\begin{equation}
\nabla_k T_i^j=\partial_k T_i^j - \Gamma^m_{ki} T_m^j + \Gamma_{mk}^j T_i^m \label {ambient C}
\end{equation}

Definition (\ref{ambient C}) provides metrilinic property of metric tensors $\nabla_\gamma X_{\alpha\beta}, \nabla_k S_{ij}=0$, consequently the surface base vectors are orthogonal to their covariant derivatives ${\bm{S}_i\cdot\nabla_k\bm{S}_j=0}$, therefore
\begin{equation}
\nabla_i\bm{S}_j=\bm{N}B_{ij} \label{curvature tensor}
\end{equation} 
$\bm{N}$ is the surface normal and $B_{ij}$ is the curvature tensor and according to definition (\ref{curvature tensor}) the curvature tensor is symmetric and extrinsic parameter. Trace of the mixed curvature tensor $B_i^i$ is the mean curvature, while the determinant is the Gaussian curvature $K=\lvert B^\cdot_\cdot \rvert$.

\section{Calculus for Moving Surfaces} All definitions from Gaussian geometry (\ref{definition} -- \ref{curvature tensor}) are also applicable to moving manifolds. However, extensions are required to parameterize the surface motion. 
The introduction of parametric time immediately adds to the necessity of defining time derivatives to preserve covariance. For that, we first introduce surface velocities: $C$ the surface normal velocity, and $V^i$ tangent ones as components of ambient velocity 
\begin{equation}
\bm{V_S}=\partial_t \bm{R}=C\bm{N}+V^i\bm{S}_i \label{surface velocity}
\end{equation}
Note here that $\bm{V_S}=(V^\alpha), C, V^i$ are only defined in ambient flat space and, therefore, are tensors only in the space. Note that the interface velocity 
\begin{equation}
C=\bm{N}\cdot\bm{V_S}=N^\alpha\bm{X}_\alpha V_\beta\bm{X}^\beta=N^\alpha V_\alpha \nonumber
\end{equation}
normal component, the surface velocity is invariant. Therefore, if the surface is assumed to be moving initially, it will retain motion in all reference frames. The invariance of interface velocity straightforwardly indicates that, in general, for moving manifolds, a reference frame cannot be found where $\bm{V_S}$ vanishes. 

Introduction of the velocity field (\ref{surface velocity}), without further elaboration (for details, see \cite{Svintradze2017, Grinfeld2013}), opens the path toward defining the covariant time derivative
\begin{equation}
\dot\nabla T_a^b=\partial_tT_a^b-V^k\nabla_k T_a^b+\dot\Gamma_k^b T^k_a-\dot\Gamma^k_a T_k^b \label{covariant t}
\end{equation}
where $\dot\Gamma_a^b=\nabla_a V^b-CB_a^b$ is the Christoffel symbols analog for invariant time derivatives first introduced by P. Grinfeld \cite{Grinfeld2013}. $\dot\nabla$ is invariant and has all properties as covariant derivatives. According to metrilinic property of covariant space/time derivatives $\nabla_k S_{ij},\dot\nabla S_{ij}=0$, therefore taking into account (\ref{covariant t}) we have
\begin{align}
&\dot\nabla S_{ij}=\partial_t S_{ij}-V^k\nabla_k S_{ij}-\dot\Gamma^k_i S_{kj}-\dot\Gamma_j^k S_{ik} \nonumber \\
&=\partial_t S_{ij}-(\nabla_i V^k-CB^k_i)S_{kj}
-(\nabla_j V^k-CB_j^k)S_{ik} \nonumber \\
&=\partial_t S_{ij}-(\nabla_i V_j-CB_{ij})-(\nabla_j V_i-CB_{ji}) \nonumber \\
&=\partial_t S_{ij}-\nabla_i V_j-\nabla_j V_i+2CB_{ij}=0, \Rightarrow \nonumber \\
&\partial_t S_{ij}=\nabla_i V_j+\nabla_j V_i-2CB_{ij} \label {metric derivative} 
\end{align}
Note that even though the invariant time derivative (\ref{covariant t}) and covariant space derivatives (\ref{ambient C}) have the same properties, they are fundamentally different because one is intrinsic. The last one is extrinsic and includes ambient surface velocities and curvature tensors.  

Definitions (\ref{ambient C}--\ref{curvature tensor}, \ref{surface velocity}, \ref{covariant t}) lead to integration theorems 
\begin {align}
\partial_t\int_\Omega Fd\Omega&=\int_\Omega\partial_t F d\Omega+\int_S FCdS \label{space integral} \\
\partial_t\int_S FdS&=\int_S \dot\nabla F dS-\int_S  FCB_i^i dS \label{surface integral}
\end{align}
Equations (\ref{definition}, \ref{covariant t}) along with integration theorems (\ref{space integral}, \ref{surface integral}) form the foundation of CMS. 
Note here that for incompressible systems, due to the conservation of volume, according to (\ref{space integral}) and Gauss theorem, one must have 
\begin{equation}
\partial_t\int_\Omega d\Omega=\int_SCdS=\int_\Omega \nabla_\alpha V^\alpha d\Omega=0 \label{incompressibility condition}
\end{equation} 
From where the condition $C,\nabla_\alpha V^\alpha=0$ follows. We frequently use (\ref{incompressibility condition}) when discussing incompressible systems. 

The proofs of (\ref{space integral}) are known in the literature. Here, we present a short and elegant proof of the (\ref{surface integral}) theorem that, to our knowledge, was not known before, though it was presented in \cite{Grinfeld2013, MGrinfeld}. 

Let $F\in S$ be a smooth scalar function on the closed smooth manifold $S$. Then, across any closed stationary contour $\gamma$, the contour integral vanishes because the contour normal $n_i$ lying in the tangent plane always stays perpendicular to the surface tangent velocity $\bm{V}_T$, so that 
\begin{equation}
v=\bm{n}\cdot\bm{V}_T=n_iV^i=0 \label{static condition}
\end{equation}
Using this and that, the total time derivative of the integral has two parts. The first part reflects that $F$ changes as time evolves, and the second part reflects that the surface changes as time evolves. Therefore 
\begin{align}
\frac{d}{dt}&\int_SFdS=\int_S\partial_tFdS+\int_SF\partial_tdS \nonumber \\
&=\int_S(\partial_tF+F\nabla_iV^i-F\nabla_iV^i)dS+\int_x\int_yF\partial_t\sqrt{S}dxdy\nonumber \\
&=\int_S(\partial_tF+F\nabla_iV^i-F\nabla_iV^i)dS\nonumber \\
&+\int_x\int_yF(\nabla_iV^i-CB_i^i)\sqrt{S}dxdy\nonumber \\
&=\int_S(\partial_tF+F\nabla_iV^i-F\nabla_iV^i)dS\nonumber \\
&+\int_SF(\nabla_iV^i-CB_i^i)dS\nonumber \\
&=\int_S(\partial_tF+\nabla_i(FV^i)-V^i\nabla_iF-CFB_i^i)dS \label{integration theorem derivation}
\end{align}
By Gauss theorem, integral from $\nabla_i(FV^i)$ is converted to the contour integral and, according to the (\ref{static condition}) boundary condition, vanishes. If one relaxes the condition that the surface is closed and is bounded by static contour, then according to (\ref{integration theorem derivation}), the integration theorem (\ref{surface integral}) updates as
\begin{equation}
\partial_t\int_SFdS=\int_S \dot\nabla F dS-\int_S  FCB_i^i dS+\int_\gamma vFd\gamma \label{contour}
\end{equation}
Since we discuss only closed surfaces here, the last term from (\ref{contour}) becomes irrelevant. 

Note that, for calculation of $\partial_t\sqrt{S}$ in (\ref{integration theorem derivation}), we used Jacobi's theorem about taking derivative for determinants ${\partial S/\partial S^{ij}=SS_{ij}}$ and (\ref{metric derivative}) we have
\begin{align}
\partial_tS=\frac{\partial S}{\partial S^{ij}}\partial_tS^{ij}&=SS_{ij}(\nabla^i V^j+\nabla^j V^i-2CB^{ij})\nonumber \\
&=2S(\nabla_iV^i-CB_i^i) \label{determinant derivative}
\end{align}
From (\ref{determinant derivative}) according to the complex function derivative rule, the identity 
\begin{equation}
\partial_t\sqrt{S}=\frac{1}{2\sqrt{S}}\partial_tS=\sqrt{S}(\nabla_iV^i-CB_i^i)
\end{equation}
trivially follows. Next, using the (\ref{covariant t}) definition written for scalar functional
 \begin{equation}
 \dot\nabla F=\partial_t F-V^i\nabla_i F \label{scalar function}
 \end{equation}
Applying (\ref{scalar function}) to the (\ref{integration theorem derivation}) with taking into account the Gauss theorem and the boundary condition (\ref{static condition}), the integration theorems (\ref{surface integral}, \ref{contour}) trivially follows.

\section{Derivation of Shape Dynamics Equations}
We start by taking the variation of the kinetic term first. In the initial steps, the derivation looks like a derivation of the NS kinetic part. We drop superscript $S$ in the index to avoid interference with tensorial indexes and recover at the end of derivations. 

To deduce the equations of motion, we first derive the simplest one from the set (\ref{MME}). According to the conservation of mass law, the following boundary conditions must be satisfied. At the end of the variation, the total mass of the surface is unchanged $dm/dt=0$. The pass integral along any closed contour $\gamma$ from the surface must vanish because contour normal velocity $v=n_iV^i=0$  ($n_i$ is the contour unit normal situated on the tangent plain) when the surface reaches static shape at the end of the variation. Considering these boundary conditions, the Gauss theorem, conservation of mass, and integration formulas (\ref{space integral}, \ref{surface integral}), we have
\begin{align}
0&=\int_\gamma\rho vd\gamma=\int_\gamma\rho n_iV^id\gamma=\int_S \nabla_i(\rho V^i)dS\nonumber \\
&=\int_S (\nabla_i(\rho V^i)+\rho CB_i^i-\rho CB_i^i)dS \nonumber \\
&=\int_S (\nabla_i(\rho V^i)-\rho CB_i^i)dS+\int_S\dot\nabla\rho dS-\frac{d}{dt}\int_S\rho dS \nonumber \\
&=\int_S(\dot\nabla\rho+\nabla_i(\rho V^i)-\rho CB_i^i)dS\label{mass conservation}
\end{align}
Since (\ref{mass conservation}) holds for every integrand, the first equation follows. To deduce the second and third equations, we note 
\begin{align}
\delta(\int_S&\Lambda_S dS+\int_\Omega P_S d\Omega)=\int_S(\dot\nabla\Lambda_S-\Lambda_SCB_i^i)dS\nonumber \\
&+\int_\Omega\partial_tP_Sd\Omega+\int_SP_SCdS \nonumber \\
&=\int_S(\partial_t\Lambda_S-V_i\nabla^i\Lambda_S-\Lambda_SCB_i^i)dS \nonumber \\
&+\int_\Omega\partial_tP_Sd\Omega+\int_SP_SCdS \label{energy variation1}
\end{align}
Next, we note that $V^i\nabla_i\Lambda_S$ governs the surface tangent deformation while all other terms in (\ref{energy variation1}) are responsible for normal deformations. Therefore, by modeling and Gauss theorem
\begin{align}
\partial_t\Lambda_S&=N_\alpha\partial_tF^\alpha \label{surface tension model} \\
\int_S\partial_t\Lambda_SdS&=\int_\Omega\partial_\alpha\partial_tF^\alpha d\Omega \label{surface tension space}
\end{align}  
Taking the (\ref{surface tension model}, \ref{surface tension space}) into account, the tangent and normal parts of the variation straightforwardly follows
\begin{align}
&\delta_T(\int_S\Lambda_S dS+\int_\Omega P_S d\Omega)=-\int_SV_i\nabla^i\Lambda_SdS\label{energy tangent} \\
&\delta_N(\int_S\Lambda_S dS+\int_\Omega P_S d\Omega)\label{energy normal} \\
&=\int_\Omega[\partial_tP_S+\partial_\alpha\partial_tF^\alpha+\partial_\alpha (P_SV^\alpha)+\partial_\alpha(\Lambda_SB_i^iV^\alpha)]d\Omega \nonumber
\end{align}
Next, let's proceed with taking the variation of the kinetic part
\begin{align}
\delta\int_S& \frac{\rho_S  V^2_S }{2}dS =\int_S(\dot\nabla\rho \cdot \frac{V^2}{2}+\dot\nabla\frac{ V^2}{2}\cdot\rho-\frac{\rho V^2}{2}CB_i^i)dS \nonumber \\
&=\int_S((\rho CB_i^i-\dot\nabla_i(\rho V^i))\cdot \frac{V^2}{2}+\dot\nabla\frac{ V^2}{2}\cdot\rho)dS \nonumber \\
&-\int_S\frac{\rho V^2}{2}CB_i^idS \nonumber \\
&=\int_S(-\nabla_i(\rho V^i \frac{V^2}{2})+\rho V^i\nabla_i\frac{V^2}{2}+\dot\nabla\frac{ V^2}{2}\cdot\rho)dS \nonumber \\
&=\int_S(\rho V^i\nabla_i\frac{V^2}{2}+\dot\nabla\frac{ V^2}{2}\cdot\rho)dS \nonumber \\ 
&=\int_S\rho \bm{V}( V^i\nabla_i \bm{V}+\dot\nabla \bm{V})dS \label{surface kinetic term}
\end{align}
Here, we must decompose (\ref{surface kinetic term}) into normal and tangent components, which requires tricky algebraic manipulations. Instead of restating the algebra here, we can refer to the direct derivation given in \cite{Svintradze2017, Svintradze2018} and give the final form only 
\begin{align}
\dot\nabla \bm{V}+V^i\nabla_i \bm{V}&=(\dot\nabla C+2V^i\nabla_iC+V^aV^bB_{ab})\bm{N}\label{dot product}\\
&+(\dot\nabla V^j+V^i\nabla_iV^j-C\nabla^jC-CV^iB_j^i)\bm{S}_j \nonumber
\end{align}
Dotting (\ref{dot product}) on the surface velocity $\bm{V}$ and taking into account the definition (\ref{surface velocity}) we get
\begin{align}
\delta\int_S &\frac{\rho_S  V^2_S }{2}dS=\int_SC(\dot\nabla C+2V^i\nabla_iC+V^aV^bB_{ab})\nonumber\\
&+V_i(\dot\nabla V^i+V^j\nabla_jV^i-C\nabla^iC-CV^jB_i^j)dS \label{integral dot product}
\end{align}
Equating (\ref{integral dot product}) with (\ref{energy variation1}), and considering the tangent and normal components of the energy variation (\ref{energy tangent}, \ref{energy normal}), and applying Gauss's theorem to the normal components, we obtain
\begin{align}
\int_\Omega\partial_\alpha[V^\alpha(&\rho_S(\dot\nabla C+2V^i\nabla_i C+V^a V^b B_{ab})\label{MME1}\\ 
&-P_S+\Lambda_S B_i^i)]d\Omega=\int_\Omega(\partial_t P_S+\partial_\alpha\partial_tF_S^\alpha)d\Omega \nonumber 
\end{align}
and
\begin{align}
\int_S\rho_SV_i(\dot\nabla V^i+V^k\nabla_k V^i-&C\nabla^i C-CV^k B_k^i)dS \nonumber \\
&=-\int_SV_i\nabla^i\Lambda_SdS\label{MME2}
\end{align}
From (\ref{mass conservation}, \ref{MME1}, \ref{MME2}) the moving manifold equations (\ref{MMEI}) trivially follow.

\section{Derivation of NS Equation}
We start by taking the variation of the kinetic term first. Because superscript $\Omega$ appears everywhere in this subsection, to minimize the number of letters in the equations and avoid interference with tensorial indexes, we drop it during calculations and recover it at the end. 

\begin{align}
\delta\int_\Omega& \frac{\rho_\Omega  V^2_\Omega }{2}d\Omega =\int_\Omega(\partial_t\rho \cdot \frac{V^2}{2}+\partial_t\frac{ V^2}{2}\cdot\rho)d\Omega \nonumber\\
&+\int_\Omega\partial_\alpha(\frac{\rho V^2}{2}V^\alpha)d\Omega=\nonumber \\
&=\int_\Omega(-\partial_\alpha(\rho V^\alpha)\cdot \frac{V^2}{2}+\partial_t\frac{ V^2}{2}\cdot\rho)d\Omega \nonumber \\
&+\int_\Omega\partial_\alpha(\frac{\rho V^2}{2}V^\alpha)d\Omega \nonumber \\
&=\int_\Omega(-\partial_\alpha(\rho V^\alpha \frac{V^2}{2})+\rho V^\alpha\partial_\alpha\frac{V^2}{2}+\partial_t\frac{ V^2}{2}\cdot\rho)d\Omega \nonumber \\
&+\int_\Omega\partial_\alpha(\frac{\rho V^2}{2}V^\alpha)d\Omega \nonumber \\
&=\int_\Omega(\rho V^\alpha\partial_\alpha\frac{V^2}{2}+\partial_t\frac{ V^2}{2}\cdot\rho)d\Omega \nonumber \\ 
&=\int_\Omega\rho_\Omega \bm{V_\Omega}( V^\alpha\partial_\alpha \bm{V_\Omega}+\partial_t \bm{V_\Omega})d\Omega \label{kinetic term1}
\end{align}
here we used the continuity equation $\partial_t\rho+\partial_\alpha(\rho V^\alpha)=0$ and $V^2=\bm{V}\cdot\bm{V}$.

Next, we model the gradient of energy, the inertial force acting on the moving fluid, therefore clarifying (\ref{momentum flux}) and its derivative. The viscous stress tensor $\sigma_{\alpha\beta}\prime$ is the only less defined term in the (\ref{momentum flux}).
\begin{equation}
\partial_\alpha E_\Omega=-\partial_\beta M_\alpha^\beta \label{energy gradient}
\end{equation}
In a fluid, internal friction processes arise when different particles of the fluid move at different velocities. This means that various parts of the fluid are in relative motion. Therefore, the viscous stress tensor depends on the velocity's space derivatives. If the velocity gradients are small, we can assume that the momentum transfer due to viscosity depends only on the first derivatives of the velocity. To the same approximation, $\sigma_{\alpha\beta}\prime$ may be a linear function of the derivatives $\partial_\alpha V_\beta$. There can be no terms $\sigma_{\alpha\beta}\prime$ independent of $\partial_\alpha V_\beta$, since the viscous stress tensor must vanish for constant velocities. Next, we notice that $\sigma_{\alpha\beta}\prime$ must also vanish when the whole fluid is in uniform rotation since it is clear that in such a motion, no internal friction occurs in the fluid. In uniform rotation with angular velocity $\omega$, the velocity equals the vector product $V=r\times\omega$. The sums are the linear combination of $\partial_\alpha V_\beta+\partial_\beta V_\alpha$ and vanish when cross-product does \cite{landau2013fluid}. Hence, the viscous stress tensor has to be proportional to linear combinations and divergence, therefore 
\begin{equation}
\sigma_{\alpha\beta}\prime=\mu(\partial_\alpha V_\beta+\partial_\beta V_\alpha-\frac{2}{3}\delta_{\alpha\beta}\partial_\gamma V^\gamma)+\xi\delta_{\alpha\beta}\partial_\gamma V^\gamma \label{stress tensor}
\end{equation}
$\mu$ is referred as viscosity and $\xi$ is called second viscosity. Taking (\ref{stress tensor}) into account in (\ref{momentum flux}) and taking the derivative of it one gets
\begin {align}
\partial_\alpha E_\Omega=-\partial_\beta M_\alpha^\beta&=-\partial_\alpha p+\partial_\beta(\mu(\partial_\alpha V^\beta+\partial^\beta V_\alpha \nonumber \\
&-\frac{2}{3}\delta_\alpha^\beta\partial_\gamma V^\gamma)+\xi\delta_\alpha^\beta\partial_\gamma V^\gamma) \label{stress tensor}
\end{align}
Taking into account that the first $\mu$ and second $\xi$ viscosities are assumed to be velocity-independent, considering that in isotropic fluids, coefficients are scalar quantities only. Note here that our calculations are model-free, meaning $\partial_\alpha E_\Omega$ can be modeled as necessary, but it does not change the structure of the kinetic terms.  Plugging (\ref{stress tensor}) in (\ref{energy gradient}) and adding (\ref{kinetic term1}) one gets 
\begin{align}
\mathcal{L}_\Omega&=\int_\Omega\rho_\Omega \bm{V_\Omega}\cdot[( V^\alpha\partial_\alpha \bm{V_\Omega}+\partial_t \bm{V_\Omega}) \nonumber \\
&-(-\bm{\nabla}p+\mu\Delta \bm{V_\Omega}+(\xi+\frac{1}{3}\mu)\bm{\nabla}\partial_\alpha V_\Omega^\alpha)]d\Omega
\end{align} 
Note that $p$ is unrelated to the surface pressure $P_S$. $p$ is the pressure induced by inertial motion, while $P_S$ is the surface pressure the system applies to the boundary, regardless of whether it is inertially moving.  

\section{Dominantly Compressible Systems}
Here, for self-consistency, we repeat the solution to (\ref{solved}), though it is published \cite{Svintradze2019}. Derivation of the condition (\ref{sphere}) $\partial_\alpha V^\alpha=0$ is trivial, and it leads to easily obtainable by algebra the solution
\begin{equation}
\bm{V}=k\frac{\bm{R}}{R^3}\Rightarrow \frac{\partial {\boldsymbol{R}}}{\partial t}={\omega }_{\xi }{R}_{\xi }{{\boldsymbol{S}}}^{{\boldsymbol{\xi }}}
 \label{spherical solution1}
\end{equation}
here $k$ is some coupling constant and $\bm{R}$ is the position vector, ${\omega }_{\xi },{R}_{\xi }$ are $\xi$ component of the frequency and position vector, ${{\boldsymbol{S}}}^{{\boldsymbol{\xi }}}$ is the $\xi$ vector component of the unit spherical vector $\bm{S}=(\sin\phi\sin\theta, \sin\phi\cos\theta, \cos\theta)$. The equation (\ref{spherical solution1}) has two trivial solutions and one, according to vector calculus, the rotor of some angular velocity like functional $\bm{L}$ solution
\begin{align}
{\boldsymbol{R}}&={{\boldsymbol{R}}}_{{\bf{0}}}+{\omega }_{\xi }{R}_{\xi }t{{\boldsymbol{S}}}^{{\boldsymbol{\xi }}}\label{S1} \\
{\boldsymbol{R}}&={R}_{0\xi }{e}^{{\omega }_{\xi }t{{\boldsymbol{S}}}^{{\boldsymbol{\xi }}}}\label{S2} \\
\bm{V}&=-\bm{\nabla}\times \bm{L} \label{rotor}
\end{align}
where $e^{\omega_\alpha \bm{S}^\alpha t}=(e^{\omega_x \bm{S}^x t},e^{\omega_y \bm{S}^y t},e^{\omega_z \bm{S}^z t})$ by definition. Next, assuming that $\bm{L}$ behaves same as the position vector (\ref{rotor}) solution can be expanded as
\begin{align}
\partial_t\bm{R}=-\bm{\nabla}\times \bm{L}\label{rotor1} \\
\frac{1}{v_0^2}\partial_t\bm{L}=\bm{\nabla}\times\bm{R}\label{rotor2}
\end{align}
$v_0$ is some constant wave propagation velocity. After some algebra on (\ref{rotor1}, \ref{rotor2}), one ends up with the wave equation for the position vector
\begin{equation}
\frac{1}{{v}_{0}^{2}}\frac{{\partial }^{2}{\boldsymbol{R}}}{\partial {t}^{2}}-{\nabla }^{2}{\boldsymbol{R}}=0 \label{waveR}
\end{equation}
The solution to the (\ref{waveR}) is straightforward, considering that all spacial directions are linearly independent, and therefore, in each direction, one has a linear wave equation. Therefore
\begin{equation}
{\boldsymbol{R}}={{\boldsymbol{R}}}_{0}+{\psi }_{\xi }{{\boldsymbol{S}}}^{{\boldsymbol{\xi }}} \label{S3}
\end{equation} 
where, the wave function ${{\boldsymbol{\psi }}}_{\xi }(\xi ,t)$ is given by the formula
\begin{align}
{{\boldsymbol{\psi }}}_{\xi }(\xi ,t)=\sum _{m=1}^{\infty }\,\frac{2}{l}({\int }_{0}^{l}\,{\boldsymbol{\psi }}(\xi ,0)\,\sin \,\frac{m\pi \xi }{l}d\xi ) \nonumber \\
\cdot\cos (\frac{{v}_{0}m\pi t}{l})\sin (\frac{{v}_{0}m\pi t}{l}) \label{S3'}
\end{align}
Some linear combination of the (\ref{S1}, \ref{S2}, \ref{S3}) solutions gives  (\ref{SF}) final solution.

To derive (\ref{wave}), we start by considering the equilibrium surface density $\rho_0$. When both $C$ and $B_i^i$ are infinitesimally small, the linearized conservation of mass can be expressed as $\partial_t\rho_S=0$, with the solution $\rho_S=\rho_0$. This suggests that for sufficiently small $C$, the density of the diffusive layer remains constant. Therefore, by taking the time derivative of the third equation of (\ref{solved}) we have
\begin{align}
\rho_0\partial_t^2C=-\Lambda_0\partial_tB_i^i&=-\Lambda_0(\nabla_i\nabla^iC+CB_{ij}B^{ij})\label{wave derivative}\\
\dot\nabla B_i^i&=\nabla_i\nabla^iC+CB_{ij}B^{ij}\label{theorem}
\end{align}
here (\ref{theorem}) is well known theorem of the CMS \cite{Grinfeld2013} and in linear approximations $\dot\nabla\approx\partial_t$. Taking into account that $B_{ij}B^{ij}=-2K$, $K$ is the Gaussian curvature \cite{Grinfeld2013}, then from (\ref{wave derivative}) directly follows 
\begin{equation}
\frac{{\partial }^{2}C}{\partial {t}^{2}}=-\gamma ({\rm{\Delta }}C-2KC)
\end{equation}
where $\gamma=\Lambda_0/\rho_0$. In planar approximations, $K\approx 0$ and the wave equation follows
\begin{equation}
\frac{{\partial }^{2}C}{\partial {t}^{2}}=-\gamma {\rm{\Delta }}C
\end{equation} 



\bibliography{apssamp}
\end{document}